\documentclass[11pt]{amsart}
\usepackage{amsmath,amssymb,amsthm,amscd, mathtools, latexsym}
\usepackage{enumerate,varioref, dsfont}
\usepackage{tikz-cd}
\usepackage{enumitem}
\usepackage[hmargin=30mm,vmargin=30mm]{geometry}
\usepackage{hyperref}
\usepackage{url}

\newtheorem{Thm}{Theorem}[section]
\newtheorem{Lem}{Lemma}[section]
\newtheorem{Prop}{Proposition}[section]

\newtheorem{Conj}{Conjecture}[section]
\newtheorem{Rem}{Remark}[section]
\newtheorem{Def}{Definition}[section]

\def\cit{{\mathbb C}}

\def\qit{{\mathbb Q}}
\def\zit{{\mathbb Z}}

\def\pit{{\mathbb P}}

\def\0{{\mathcal O}}

\def\Hom{\mathop{\rm Hom}\nolimits}

\def\h{{\mathfrak h}}

\def\h{{\mathfrak h}}

\def\M{{\mathcal M}}

\def\V{{\mathcal V}}

\begin{document}

\title{THE MOTIVE OF A SMOOTH THETA DIVISOR}

\author{Humberto A. Diaz} 
\thanks{The author would like to thank Chad Schoen, his advisor, and Bruno Kahn for their help with edits.}

\address{DEPARTMENT OF MATHEMATICS, DUKE UNIVERSITY, DURHAM, NC}
\email{hdiaz123@math.duke.edu}

\date{}

\maketitle

\begin{abstract}
We prove a motivic version of the Lefschetz hyperplane theorem for a smooth ample divisor $\Theta$ on an Abelian variety. We use this to construct a motive $P$ that realizes the primitive cohomology of $\Theta$.  \end{abstract}

\section{Introduction}

Let $k$ be an algebraically closed field. Given a smooth projective variety $X$ of dimension $d$ over $k$ and a Weil cohomology $H^{*}$, there is a decomposition of the diagonal $[\Delta_{X}] \in H^{2d} (X \times X)$ into its K{\"u}nneth components:
\begin{equation} \Delta_{X} = \Delta_{0, X} + ...+ \Delta_{2d, X} \in H^{2d} (X \times X) \cong \bigoplus H^{j} (X) \otimes H^{2d-j} (X) \nonumber \end{equation}

\noindent It is one of Grothendieck's standard conjectures (\cite{KL} Section 4) that these K{\"u}nneth components arise from algebraic cycles; i.e., that there exist correspondences $\pi_{j, X} \in CH^{d} (X \times X)$ for which $cl(\pi_{j, X}) = \Delta_{j, X}$ under the cycle class map $cl: CH^{d} (X \times X) \to H^{2d} (X \times X)$. We can state a stronger version of this conjecture as follows:
\begin{Conj}[Chow-K{\"u}nneth]\label{CK} There exist correspondences $\pi_{j, X} \in CH^{d} (X \times X)$ satisfying:
\begin{enumerate}[label=(\alph*)] \item\label{mut-orth} $\pi_{j, X}^{2} = \pi_{j, X}, \ \pi_{j, X}\circ \pi_{j', X}  = 0$ for $j \neq j'$  \item\label{sum} $\sum \pi_{j, X} = \Delta_X$ \item\label{Weil} $cl(\pi_{j, X}) = \Delta_{j, X}$ {for any choice of Weil cohomology.} \end{enumerate} \end{Conj}

\noindent In this stronger version, the correspondences $\pi_{j, X}$ are actually idempotents, which gives Chow motives $\h^{j} (X) = (X, \pi_{j, X}, 0)$. Moreover, the decomposition of the diagonal into {\em orthogonal components} gives a decomposition of the motive of $X$ as $\bigoplus\h^{j} (X)$. An important problem in the theory of motives is to understand these ``underlying'' objects $\h^{j} (X)$ that represent the various degrees of cohomology (for every choice of cohomology). The Chow-K{\"u}nneth conjecture is known to hold in some important cases: curves, surfaces (\cite{M}), complete intersections in $\pit^{n}$ (\cite{MNP} Chapter 6), Abelian varieties (\cite{DM}), elliptic modular varieties (\cite{GHM}).

\indent Suppose that $A$ is an Abelian variety of dimension $g$ and $i: \Theta \xhookrightarrow{} A$ is a smooth ample divisor. The first goal of this note is then to prove the following:

\begin{Thm}\label{1} There exist correspondences $\pi_{j, \Theta} \in CH^{g-1} (\Theta \times \Theta)$ satisfying conjecture \ref{CK}. 
\end{Thm}

The Lefschetz hyperplane theorem gives isomorphisms $i^{*} : H^{j} (A) \to H^{j} (\Theta)$ for $j < g-1$ and $i_{*} : H^{j} (\Theta) \to H^{j+2} (A)$ for $j > g-1$. The proof of Theorem \ref{1} gives a particular set of idempotents $\pi_{j, \Theta}$ and we set $\h^{j} (\Theta) = (\Theta, \pi_{j, \Theta}, 0)$. We also set $\h^{j} (A) = (A, \pi_{j, A}, 0)$, where $\pi_{j, A}$ are the canonical idempotents constructed in \cite{DM}. We are then able to prove the following motivic version of the Lefschetz hyperplane theorem: 
\begin{Thm}\label{2}
\begin{enumerate}[label=(\alph*)] \item\label{j<} The pull-back $\h^{j} (i) := \pi_{j, \Theta} \circ \prescript{t}{}{\Gamma_{i}} \circ \pi_{j, A}: \h^{j} (A) \rightarrow \h^{j} (\Theta)$ is an isomorphism for $j < g-1$.
\item\label{j>} The push-forward $\prescript{t}{}{\h^{j} (i)} := \pi_{j+2, A} \circ \Gamma_{i} \circ \pi_{j, \Theta}: \h^{j} (\Theta)  \rightarrow \h^{j+2} (A)(1)$ is an isomorphism for $j > g-1$.
\item\label{j=} $\h^{g-1} (i)$ is split-injective and $\prescript{t}{}{\h^{g-1} (i)}$ is split-surjective.
\item\label{P} There is an idempotent $p \in CH^{g-1} (\Theta \times \Theta)$ which is orthogonal to $\pi_{j, \Theta}$ for $j \neq g-1$ and for which the motive $P := (\Theta, p, 0)$ satisfies $H^{*} (P) = K_{\Theta} := \ker(i_{*}: H^{g-1} (\Theta) \to H^{g+1} (A))$.
\end{enumerate}
\end{Thm}

We can specialize to the case that $k = \cit$ and $H^{*}$ is singular cohomology with $\qit$-coefficients. The primitive cohomology of $\Theta$, $$K_{\Theta} = \ker(i_{*}: H^{g-1} (\Theta, \qit) \to H^{g+1} (A, \qit) (1)),$$ is the only Hodge substructure of $H^{*}$ not coming from $A$. So, one should expect to encounter difficulty in analyzing the motive $P$. The simplest nontrivial case is when $A$ is a principally polarized Abelian fourfold and $[\Theta] \in CH^{1} (A)$ is its principal polarization. In this case, $\Theta$ is generally a smooth divisor and $H^*(P) = K_{\Theta}$ has Hodge level 1. Conjecturally, a motive over $\cit$ whose singular cohomology has Hodge level 1 should correspond to an Abelian variety (\cite{KI} Remark 7.12). We have the following partial result:
\begin{Prop}\label{3} There exists an Abelian variety $J$ such that $\h^{1} (J) (-1) \cong P$ $\Leftrightarrow p_{*}CH_{0} (\Theta) = 0$.
\end{Prop}
 
\section{Preliminaries}
Let $\M_{k}$ denote the category of Chow motives over $k$ whose objects are triples $(X, \pi, n)$, where $X$ is a smooth projective variety of dimension $d$, $\pi \in CH^{d} (X \times X)$ is an idempotent and $n \in \zit$. The morphisms are defined as follows:
\begin{equation} \begin{split}
\Hom_{\M_k} ((X, \pi, n), (X', \pi', n')) & := \pi'\circ Cor^{n'-n} (X, X')\circ \pi \\ 
&= \pi'\circ CH^{d+ n'-n} (X \times X')\circ \pi \nonumber \end{split}
\end{equation}
Here, composition is defined in \cite{F} Chapter 16.1. There is a functor $\mathfrak{h}: \V_k^{opp} \mapsto \M_k$ from the category of smooth projective varieties over $k$ with $\mathfrak{h} (X) = (X, \Delta_{X}, 0)$ and with $\mathfrak{h} (g) = \prescript{t}{}{\Gamma_{g}}$ for any morphism $g: X \to X'$. A Weil cohomology theory is a functor $H^{*} : \V_{k}^{opp} \mapsto Vec_{K}$ (with $K$ is a field of characteristic 0) satisfying certain axioms (described in \cite{KL} Section 4), one of which is the Lefschetz hyperplane isomorphism. Examples include singular, $\ell$-adic, crystalline, or de Rham cohomology. This extends to a functor $H^{*} : \M_{k} \mapsto Vec_{K}$, and for $M = (X, \pi, m)$, we have $$H^{j} (M) = \pi_{*}H^{j+2m} (X).$$Also, there is the extension of scalars functor $( )_{L}: \M_{k} \mapsto \M_{L}$ for any field extension $k \subset L$. For $M = (X, \pi, m)$, we will use the notation $M(n) = (X, \pi, m-n)$ and $\mathds{1} (n) = (\text{Spec } k, \Delta_k, n)$. 

\begin{Lem}[Liebermann]\label{Lieb} Let $h_{X}: X' \vdash X$, $h_{Y}: Y \vdash Y'$ be correspondences of smooth projective varieties. Then, for $\alpha \in CH^{*} (X \times Y), \beta \in CH^{*} (X' \times Y')$, we have 
\begin{enumerate}[label=(\alph*)]
\item\label{itm: h} $(h_{X} \times h_{Y})_{*} (\alpha) = h_{Y} \circ \alpha \circ h_{X}$
\item\label{itm: a} When $f: X \to X'$ and $g: Y \to Y'$ are morphisms, $(f \times g)_{*} (\alpha) = {\Gamma_{g}} \circ \alpha \circ \prescript{t}{}{\Gamma_f}.$
\end{enumerate}
\begin{proof} See \cite{F} Proposition 16.1.1. 
\end{proof}
\end{Lem}

\begin{Thm}[Shermenev, Deninger-Murre] \label{DM} Let $A$ be an Abelian variety of dimension $g$ over $k$. Then, there is a unique set of idempotents $\{ \pi_{j, A} \} \in  CH^{g} (A \times A)$ satisfying conjecture \ref{CK} and the following relation for all $n \in \zit$:
\begin{equation} \prescript{t}{}{\Gamma_{n}}\circ\pi_{j, A} = n^j\cdot \pi_{j, A} = \pi_{j, A}\circ\prescript{t}{}{\Gamma_{n}} \label{n} \end{equation}  
\begin{proof} See \cite{DM} Theorem 3.1. \end{proof}
\end{Thm}

Let $i: \Theta \hookrightarrow A$ be a smooth ample divisor and let $\h^{j} (A) = (A, \pi_{j, A}, 0)$ be the motive for the idempotents in Theorem \ref{DM}. Then, we define the {\em Lefschetz operator}: $$L_{\Theta} := \Delta_{*}(\Theta) \in CH^{g+1} ({A} \times {A}).$$ 
The most essential result for the proofs of theorems \ref{1} and \ref{2} is the following in \cite{K}, a motivic version of the Hard Lefschetz theorem:
\begin{Thm}[K{\"u}nnemann]\label{K} Assume that $[\Theta] = (-1)_{A}^*[\Theta] \in CH^{1} (A)$.
\begin{enumerate}[label=(\alph*)]
\item\label{itm: 0} $(L_{\Theta})_{*}\alpha = \alpha\cup[\Theta]$ for $\alpha \in H^{*} (A)$
\item\label{itm: hard-lef} The operator $\pi_{2g-j, A} \circ L_{\Theta}^{g-j} \circ \pi_{j, A}: \h^{j} (A) (g-j) \to \h^{2g-j} (A)$ is an isomorphism of motives for $j \leq g$. That is, there exists a correspondence $\Lambda_{\Theta} \in CH^{g-1} ({A} \times {A})$ such that the following relations hold for $j \leq g$:
\begin{equation} \begin{split}\pi_{j, A}\circ\Lambda_{\Theta}^{g-j} \circ L_{\Theta}^{g-j} \circ \pi_{j, A} & = \pi_{j, A} \\ \pi_{2g-j, A}\circ L_{\Theta}^{g-j} \circ\Lambda_{\Theta}^{g-j} \circ \pi_{2g-j, A} & = \pi_{2g-j, A} \end{split} \label{crucial} \end{equation}
\item\label{itm: com} Set $\pi_{j, A} = 0$ for all $j \notin \{0, 1, ...2g\}$. Then, we have $L_{\Theta}\circ\pi_{j, A} = \pi_{j+2, A}\circ L_{\Theta}$ and ${\Lambda}_{\Theta}\circ\pi_{j, A} = \pi_{j-2, A}\circ {\Lambda}_{\Theta}$.
\end{enumerate}
\begin{proof}
See \cite{K} Theorem 4.1. It should be noted that \ref{itm: hard-lef} holds more generally for Abelian schemes. It is a technical result that uses properties of the Fourier transform for Chow groups of Abelian schemes. 
\end{proof} 
\end{Thm}

By Theorem \ref{K} \ref{itm: 0} and the projection formula, we have $(L_{\Theta})_{*} = \cup[\Theta] = i_{*}\circ i^{*}$. The result below shows that this is true on the level of correspondences:
\begin{Lem}\label{cor-lieb} $L_{\Theta} = \Gamma_{i}\circ\prescript{t}{}{\Gamma_{i}} \in CH^{g+1} (A \times A).$
\begin{proof} From the obvious commutative diagram:
\begin{equation}\begin{CD}
\Theta @>\Delta_{\Theta}>> \Theta \times \Theta\\
@ViVV @Vi \times iVV\\
A @>\Delta_{A}>> A \times A
\end{CD}\end{equation}
we have $L_{\Theta} = (\Delta_{A})_{*}(\Theta) = (\Delta_{A})_{*}(i_{*}1) = (i \times i)_{*}(\Delta_{\Theta}) = \Gamma_{i} \circ \Delta_{\Theta} \circ \prescript{t}{}{\Gamma_{i}} = \Gamma_{i}\circ\prescript{t}{}{\Gamma_{i}}$, where the penultimate step follows from Lemma \ref{Lieb} \ref{itm: a}. \end{proof} \end{Lem}

\section{Proofs of Theorems \ref{1} and \ref{2}}

Since $k$ is algebraically closed, it's possible to find some $a \in A(k)$ such that $t_{a}^{*}[\Theta] \in CH^{1} (A)$ is invariant under $(-1)^*$. So, we can assume that $(-1)_{A}^*[\Theta] = [\Theta]$, so that the results of the previous section are applicable.

\begin{proof}[Proof of Theorem \ref{1}] For the proof, we will need to exhibit correspondences $\pi_{j, \Theta} \in CH^{g-1} (\Theta \times \Theta)$ which satisfy conjecture \ref{CK}. These are given as follows:
\begin{equation} \begin{split} \pi_{j, \Theta} & = \prescript{t}{}{\Gamma_i}\circ\pi_{j, A}\circ\Lambda_{\Theta}^{g-j}\circ L_{\Theta}^{g-j-1}\circ\Gamma_i \text{ for } j< g-1,\\
\pi_{j, \Theta} & = \prescript{t}{}{\Gamma_i}\circ L_{\Theta}^{j-g+1}\circ \Lambda_{\Theta}^{j-g+2}\circ\pi_{j+2, A}\circ\Gamma_i \text{ for } j> g-1,\\
\displaystyle \pi_{g-1, \Theta} & = \Delta_{\Theta} - \sum_{j \neq g-1} \pi_{j, \Theta}.
\end{split} \end{equation}

Since $\sum \pi_{j, \Theta} = \Delta_{\Theta}$ holds by definition, it suffices to check conditions \ref{mut-orth} and \ref{Weil} of conjecture \ref{CK}. For $j < g-1$, we have
\begin{equation}\begin{split}
 \pi_{j, \Theta}^{2} & = \prescript{t}{}{\Gamma_i}\circ\pi_{j, A}\circ\Lambda_{\Theta}^{g-j}\circ L_{\Theta}^{g-j-1}\circ\Gamma_i \circ \prescript{t}{}{\Gamma_i}\circ\pi_{j, A}\circ\Lambda_{\Theta}^{g-j}\circ L_{\Theta}^{g-j-1}\circ\Gamma_i\\
 & = \prescript{t}{}{\Gamma_i}\circ\pi_{j, A}\circ\Lambda_{\Theta}^{g-j}\circ L_{\Theta}^{g-j}\circ\pi_{j, A}\circ\Lambda_{\Theta}^{g-j}\circ L_{\Theta}^{g-j-1}\circ\Gamma_i\\
  & = \prescript{t}{}{\Gamma_i}\circ\pi_{j, A}\circ\Lambda_{\Theta}^{g-j}\circ L_{\Theta}^{g-j-1}\circ\Gamma_i =  \pi_{j, \Theta}\nonumber
 \end{split}\end{equation}
 Here, the second equality holds by Lemma \ref{cor-lieb}, the third holds by Theorem \ref{K} \ref{itm: hard-lef}. Similarly, for $j > g-1$ we have:
 \begin{equation}\begin{split}
 \pi_{j, \Theta}^{2} & = \prescript{t}{}{\Gamma_i}\circ\pi_{j, A}\circ L_{\Theta}^{j-g+1}\circ \Lambda_{\Theta}^{j-g+2}\circ\pi_{j+2, A}\circ\Gamma_i\circ\prescript{t}{}{\Gamma_i}\circ L_{\Theta}^{j-g+1}\circ \Lambda_{\Theta}^{j-g+2}\circ\pi_{j+2, A}\circ\Gamma_i\\
 & =\prescript{t}{}{\Gamma_i}\circ L_{\Theta}^{j-g+1}\circ \Lambda_{\Theta}^{j-g+2}\circ\pi_{j+2, A}\circ L_{\Theta}\circ L_{\Theta}^{j-g+1}\circ \Lambda_{\Theta}^{j-g+2}\circ\pi_{j+2, A}\circ\Gamma_i\\
 &  =\prescript{t}{}{\Gamma_i}\circ L_{\Theta}^{j-g+1}\circ \Lambda_{\Theta}^{j-g+2}\circ \pi_{j+2, A}\circ L_{\Theta}^{j-g+2}\circ \Lambda_{\Theta}^{j-g+2}\circ \pi_{j+2, A}\circ\Gamma_i\\
  &  =\prescript{t}{}{\Gamma_i}\circ L_{\Theta}^{j-g+1}\circ \Lambda_{\Theta}^{j-g+2}\circ \pi_{j+2, A}\circ\Gamma_i =  \pi_{j, \Theta}
 \nonumber\end{split}\end{equation}
Thus, $\pi_{j, \Theta}^{2} = \pi_{j, \Theta}$ for $j \neq g-1$. Before proving the same for $j = g-1$, we show that the orthogonality condition of \ref{mut-orth} (in conjecture \ref{CK}) holds; that is, $\pi_{j, \Theta}\circ \pi_{j', \Theta} = 0$ for $j \neq j'$ and $j, j' \neq g-1$. We do this for the case of $j \neq j' < g-1$:
 \begin{equation}\begin{split}
 \pi_{j, \Theta}\circ \pi_{j', \Theta} & = \prescript{t}{}{\Gamma_i}\circ\pi_{j, A}\circ\Lambda_{\Theta}^{g-j}\circ L_{\Theta}^{g-j-1}\circ\Gamma_i \circ \prescript{t}{}{\Gamma_i}\circ\pi_{j', A}\circ\Lambda_{\Theta}^{g-j'}\circ L_{\Theta}^{g-j'-1}\circ\Gamma_i\\
 & = \prescript{t}{}{\Gamma_i}\circ\pi_{j, A}\circ\Lambda_{\Theta}^{g-j}\circ L_{\Theta}^{g-j}\circ\pi_{j', A}\circ\Lambda_{\Theta}^{g-j'}\circ L_{\Theta}^{g-j'-1}\circ\Gamma_i\\
 & = \prescript{t}{}{\Gamma_i}\circ\Lambda_{\Theta}^{g-j}\circ L_{\Theta}^{g-j}\circ\pi_{j, A}\circ\pi_{j', A}\circ\Lambda_{\Theta}^{g-j'}\circ L_{\Theta}^{g-j'-1}\circ\Gamma_i\\
  & = 0
\label{similar} \end{split}\end{equation}
Again, the second equality holds by Lemma \ref{cor-lieb} and the last equality follows from the orthogonality condition in Theorem \ref{DM}. The third equality holds because we have \begin{equation} \pi_{j, A}\circ\Lambda_{\Theta}^{g-j}\circ L_{\Theta}^{g-j} = \Lambda_{\Theta}^{g-j}\circ L_{\Theta}^{g-j}\circ\pi_{j, A}\nonumber\end{equation} which follows by repeated application of Theorem \ref{K} \ref{itm: com}. The remaining cases of orthogonality ($j \neq j'$ and $j, j' \neq g-1$) are identical to (\ref{similar}).\\
\indent What remains for the verification of condition \ref{mut-orth} is to show that:
\begin{enumerate}[label=(\roman*)] \item\label{first} $\pi_{g-1, \Theta}^{2} = \pi_{g-1, \Theta}$ \item\label{second} $\pi_{g-1, \Theta}\circ\pi_{j, \Theta} = 0 = \pi_{j, \Theta}\circ\pi_{g-1, \Theta}$ for $j\neq g-1$ \end{enumerate} For \ref{first}, let $\pi = \sum_{k \neq g-1} \pi_{j, \Theta}$. Since the summands are mutually orthogonal idempotents by the preceding verifications, it follows that $\pi^{2} = \pi$. Since $\pi_{g-1, \Theta} = \Delta_{\Theta} - \pi$ by definition, we have \begin{equation} \displaystyle \pi_{g-1, \Theta}^{2}  = (\Delta_{\Theta} - \pi)^{2} = \Delta_{\Theta} + \pi^{2} - 2\pi  = \Delta_{\Theta} - \pi = \pi_{g-1, \Theta}\nonumber\end{equation}
For \ref{second}, let $j\neq g-1$ and note that \begin{equation} \begin{split} \pi_{g-1, \Theta}\circ\pi_{j, \Theta} = (\Delta_{\Theta} - \pi)\circ \pi_{j, \Theta} & = \pi_{j, \Theta} - \sum_{k \neq g-1} \pi_{k, \Theta}\circ \pi_{j, \Theta} \\ & = \pi_{j, \Theta} - \pi_{j, \Theta} = 0\nonumber \end{split}\end{equation}
where the third equality holds since $\pi_{k, \Theta}\circ \pi_{j, \Theta} = 0$ for $j \neq k$. Similarly, one has $0 = \pi_{j, \Theta}\circ\pi_{g-1, \Theta}$. This completes the verification of item \ref{mut-orth} in conjecture \ref{CK}.\\
\indent Finally, we prove \ref{Weil} in conjecture \ref{CK}. It suffices to show that $\pi_{j, \Theta}$ acts as the identity on $H^{j} (\Theta)$ and trivially on $H^{j'} (\Theta)$ for $j \neq j'$ and any Weil cohomology $H^{*}$. One easily reduces this to the case that $j \neq g-1$. We will verify this for $j < g-1$. Since $\pi_{j, A}$ acts as 0 on $H^{j'} (A)$ for $j \neq j'$, we need only show that $\pi_{j, \Theta}$ acts as the identity on $H^{j} (\Theta)$. To this end, let $\phi := \Lambda_{\Theta}^{g-j}\circ L_{\Theta}^{g-j-1}\circ\Gamma_i$ so that $$\pi_{j, \Theta} = \prescript{t}{}{\Gamma_i}\circ\pi_{j, A}\circ\phi$$ Since $H^{*}$ is a Weil cohomology, $\prescript{t}{}{\Gamma_i}_{*} = i^{*}: H^{j} (A) \to H^{j} (\Theta)$ is an isomorphism (see \cite{KL}). Moreover, by Hard Lefschetz, $(\phi\circ\prescript{t}{}{\Gamma_i})_{*} = (\Lambda_{\Theta}^{g-j})_{*}\circ (L_{\Theta}^{g-j})_{*}$ is the identity on $H^{j} (A)$. Thus, $i^{*}$ and $\phi_{*}$ are inverses, from which it follows that $(\pi_{j, \Theta})_{*}$ is the identity on $H^{j} (\Theta)$ for $j < g-1$. The case of $j > g-1$ is nearly identical, only that one uses the fact that $i_{*}$ is an isomorphism. \end{proof}

\begin{proof}[Proof of Theorem \ref{2}] The statements of \ref{j<} and \ref{j>} are that $\h^{j} (i)$ and $\prescript{t}{}{\h^{j} (i)}$ are isomorphisms for $j < g-1$ and $j > g-1$, respectively. To show this, we need to construct their inverse isomorphisms:
\begin{equation}\begin{split}
\phi_{j} & := \pi_{j, A}\circ\Lambda_{\Theta}^{g-j}\circ L_{\Theta}^{g-j-1}\circ \Gamma_{i}\circ\pi_{j, \Theta} \text{ for } j < g-1\\
\phi_{j} & :=  \pi_{j,\Theta}\circ\prescript{t}{}{\Gamma_i}\circ L_{\Theta}^{j-g+1}\circ\Lambda_{\Theta}^{j-g+2}\circ\pi_{j+2, A}\text{ for } j > g-1
\nonumber\end{split}\end{equation}
Then, for $j < g-1$, we have 
\begin{equation}\begin{split}
 \phi_{j}\circ\h^{j} (i) & = \pi_{j, A}\circ\Lambda_{\Theta}^{g-j}\circ L_{\Theta}^{g-j-1}\circ \Gamma_{i}\circ\pi_{j, \Theta}\circ\prescript{t}{}{\Gamma_i}\circ\pi_{j, A}\\
 & = \pi_{j, A}\circ\Lambda_{\Theta}^{g-j}\circ L_{\Theta}^{g-j-1}\circ \Gamma_{i}\circ\prescript{t}{}{\Gamma_i}\circ\pi_{j, A}\circ\Lambda_{\Theta}^{g-j}\circ L_{\Theta}^{g-j-1}\circ \Gamma_{i}\circ\prescript{t}{}{\Gamma_i}\circ\pi_{j, A}\\
& = \pi_{j, A}\circ\Lambda_{\Theta}^{g-j}\circ L_{\Theta}^{g-j}\circ \pi_{j, A}\circ\Lambda_{\Theta}^{g-j}\circ L_{\Theta}^{g-j}\circ \pi_{j, A}\\
& = \pi_{j, A}
\nonumber\end{split}\end{equation}
where the third and fourth equalities hold by Theorem \ref{K} \ref{itm: hard-lef}. Similarly, we have
\begin{equation}\begin{split}
 \h^{j} (i)\circ\phi_{j} & = \pi_{j, \Theta}\circ\prescript{t}{}{\Gamma_i}\circ\pi_{j, A}\circ\Lambda_{\Theta}^{g-j}\circ L_{\Theta}^{g-j-1}\circ \Gamma_{i}\circ\pi_{j, \Theta}\\
 & = \pi_{j, \Theta}^{3} = \pi_{j, \Theta}
\nonumber\end{split}\end{equation}
We conclude that $\h^{j} (i)$ and $\phi_{j}$ are inverses for $j < g-1$, proving \ref{j<}. For \ref{j>}, we have 
\begin{equation}\begin{split}
 \prescript{t}{}{\h^{j} (i)}\circ\phi_{j} & =  \pi_{j+2, A}\circ\Gamma_i\circ\pi_{j,\Theta}\circ\prescript{t}{}{\Gamma_i}\circ L_{\Theta}^{j-g+1}\circ\Lambda_{\Theta}^{j-g+2}\circ\pi_{j+2, A}\\
& =  \pi_{j+2, A}\circ\Gamma_i\circ\prescript{t}{}{\Gamma_i}\circ L_{\Theta}^{j-g+1}\circ\Lambda_{\Theta}^{j-g+2}\circ\pi_{j+2, A}\circ\Gamma_i\circ\prescript{t}{}{\Gamma_i}\circ L_{\Theta}^{j-g+1}\circ\Lambda_{\Theta}^{j-g+2}\circ\pi_{j+2, A}\\
& =  \pi_{j+2, A}\circ L_{\Theta}^{j-g+2}\circ\Lambda_{\Theta}^{j-g+2}\circ\pi_{j+2, A}\circ L_{\Theta}^{j-g+2}\circ\Lambda_{\Theta}^{j-g+2}\circ\pi_{j+2, A}\\
& = \pi_{j, A}\nonumber
\nonumber\end{split}\end{equation}
Similarly, we have 
\begin{equation}\begin{split}
\phi_{j}\circ \prescript{t}{}{\h^{j} (i)} & = \pi_{j,\Theta}\circ\prescript{t}{}{\Gamma_i}\circ L_{\Theta}^{j-g+1}\circ\Lambda_{\Theta}^{j-g+2}\circ\pi_{j+2, A}\circ\Gamma_i\circ\pi_{j,\Theta}\\
& = \pi_{j,\Theta}^{3} =\pi_{j,\Theta}
\nonumber\end{split}\end{equation}
So, $\prescript{t}{}{\h^{j} (i)}$ and $\phi_{j}$ are inverses for $j > g-1$. For \ref{j=}, we need to show that $\h^{g-1} (i)$ and $\prescript{t}{}{\h^{g-1} (i)}$ are split-injective and split-surjective, respectively. Their left and right inverses will be:
\begin{equation}\begin{split}
\phi_{g-1} & = \pi_{g-1, A}\circ\Lambda_{\Theta}\circ\Gamma_{i}\circ\pi_{g-1, \Theta}\\ 
\psi_{g-1} & = \pi_{g-1, \Theta}\circ\prescript{t}{}{\Gamma_{i}}\circ\Lambda_{\Theta}\circ\pi_{g+1, A}.
\end{split}\end{equation}
To this end, we begin by noting that for $j < g-1$:
\begin{equation}\begin{split}
\pi_{j, \Theta}\circ\prescript{t}{}{\Gamma_{i}} & = \prescript{t}{}{\Gamma_i}\circ\pi_{j, A}\circ\Lambda_{\Theta}^{j-g}\circ L_{\Theta}^{g-j} \\
& = \prescript{t}{}{\Gamma_{i}}\circ\pi_{j, A}\label{<}
\end{split}\end{equation}
Similarly, we have ${\Gamma_{i}}\circ\pi_{j, \Theta} = \pi_{j+2, A}\circ{\Gamma_{i}}$ for $j > g-1$. So, we write $\pi = \displaystyle \sum_{j \neq g-1} \pi_{j, \Theta}$ as before and obtain:
\begin{equation} \begin{split} \Gamma_{i}\circ\pi\circ\prescript{t}{}{\Gamma_i}\circ\pi_{g-1, A} & = \sum_{j < g-1} \Gamma_{i}\circ\pi_{j, \Theta}\circ\prescript{t}{}{\Gamma_i}\circ\pi_{g-1, A} + \sum_{j > g-1} \Gamma_{i}\circ\pi_{j, \Theta}\circ\prescript{t}{}{\Gamma_i}\circ\pi_{g-1, A}\\
& = \sum_{j < g-1} \Gamma_{i}\circ\prescript{t}{}{\Gamma_i}\circ\pi_{j, A}\circ\pi_{g-1, A} + \sum_{j > g-1} \pi_{j+2, A}\circ\Gamma_{i}\circ\prescript{t}{}{\Gamma_i}\circ\pi_{g-1, A}\\
& = \sum_{j < g-1} L_{\Theta}\circ\pi_{j, A}\circ\pi_{g-1, A} + \sum_{j > g-1} \pi_{j+2, A}\circ L_{\Theta}\circ\pi_{g-1, A}\\
& = \sum_{j > g-1} L_{\Theta}\circ\pi_{j, A}\circ \pi_{g-1, A} = 0
\label{red}\end{split} \end{equation} 
where the third equality holds by the mutual orthogonality of $\pi_{j, A}$ and the fourth holds because $L_{\Theta}\circ\pi_{j, A} = \pi_{j+2, A}\circ L_{\Theta}$. Thus, we have:
\begin{equation}\begin{split}
\phi_{g-1}\circ\h^{g-1}(i) & = \pi_{g-1, A}\circ\Lambda_{\Theta}\circ\Gamma_{i}\circ\pi_{g-1, \Theta}\circ\prescript{t}{}{\Gamma_i}\circ\pi_{g-1, A}\\
& = \pi_{g-1, A}\circ\Lambda_{\Theta}\circ\Gamma_{i}\circ(\Delta_{\Theta} - \pi)\circ\prescript{t}{}{\Gamma_i}\circ\pi_{g-1, A} \\ 
& = \pi_{g-1, A}\circ\Lambda_{\Theta}\circ\Gamma_{i}\circ\prescript{t}{}{\Gamma_i}\circ\pi_{g-1, A} - \pi_{g-1, A}\circ\Lambda_{\Theta}\circ\Gamma_{i}\circ \pi\circ\prescript{t}{}{\Gamma_i}\circ\pi_{g-1, A} \\
& = \pi_{g-1, A}\circ\Lambda_{\Theta}\circ L_{\Theta}\circ\pi_{g-1, A} = \pi_{g-1, A}
\nonumber
\end{split}\end{equation}
Here, the second term on the third line vanishes by (\ref{red}). So, $\h^{g-1}(i)$ is split-injective. A similar calculation shows that $\prescript{t}{}{\h^{g-1}(i)}$ is split-surjective with right inverse $\psi_{j}$. The completes the proof of \ref{j=}. \\
\indent Finally, for \ref{P} we define:
\begin{equation}
\pi_{g-1, \Theta}' := \prescript{t}{}{\Gamma_i}\circ\pi_{g-1, A}\circ\Lambda_{\Theta}\circ\Gamma_i \in CH^{g-1} (\Theta \times \Theta) \nonumber
\end{equation}
As in the proof of Theorem \ref{1}, one can show that $\pi_{g-1, \Theta}'$ is an idempotent, is orthogonal to $\pi_{j, \Theta}$ for $j \neq g-1$. It follows that 
\begin{equation}
\pi_{g-1, \Theta}'\circ\pi_{g-1, \Theta} = \pi_{g-1, \Theta}' - \sum_{j \neq g-1} \pi_{g-1, \Theta}'\circ\pi_{j, \Theta} =\pi_{g-1, \Theta}' \nonumber
\end{equation}
Similarly, one has $\pi_{g-1, \Theta}\circ\pi_{g-1, \Theta}' = \pi_{g-1, \Theta}'$. Write $\h^{g-1}_{1} (\Theta)= (\Theta, \pi_{g-1, \Theta}', 0)$ for the corresponding motive and define:
\begin{equation}
p := \pi_{g-1, \Theta} - \pi_{g-1, \Theta}' \in CH^{g-1} (\Theta \times \Theta) \nonumber
\end{equation}
We have \begin{equation} \begin{split} p^{2} = (\pi_{g-1, \Theta} - \pi_{g-1, \Theta}')^{2} & = \pi_{g-1, \Theta}^{2} + (\pi_{g-1, \Theta}')^{2} - 2\pi_{g-1, \Theta}\circ\pi_{g-1, \Theta}'\\ & = \pi_{g-1, \Theta} + \pi_{g-1, \Theta}' - 2\pi_{g-1, \Theta}' = \pi_{g-1, \Theta} - \pi_{g-1, \Theta}' = p \nonumber\end{split}\end{equation} so that $p$ is an idempotent. Write $P:=  (\Theta, p, 0)$ for the corresponding motive. We also have $$p \circ \pi_{g-1, \Theta}' = (\pi_{g-1, \Theta} - \pi_{g-1, \Theta}')\circ \pi_{g-1, \Theta}' = \pi_{g-1, \Theta}' - \pi_{g-1, \Theta}' = 0$$ so that $p$ and $\pi_{g-1, \Theta}'$ are orthogonal. This gives a decomposition of motives: \begin{equation} \h^{g-1} (\Theta) = P \oplus  \h^{g-1}_{1} (\Theta) \label{last} \end{equation}
The same argument for Theorem \ref{1} \ref{Weil} shows that $H^{*} (\h^{g-1}_{1} (\Theta)) = i^{*}H^{g-1} (\Theta)$. Thus, applying $H^{*}$ to (\ref{last}), it follows that $H^{*} (P) = K_{\Theta}$. 
\end{proof}
\section{The complementary motive P}
\indent Now, let $k = \cit$ and $H^{*}$ be singular cohomology with $\qit$-coefficients. We consider the case of $A$ a prinicipally polarized Abelian variety, whose principal polarization is the class of $i: \Theta \to A$. Since we are interested in the motive $P$, we need $\Theta$ to be nonsingular. The simplest nontrivial case is that of $g=4$, where a well-known result of Mumford in \cite{M2} is that $\Theta$ is generally nonsingular. Now, let $K_{\Theta, \qit} := \ker(i_{*}: H^{g-1} (\Theta, \qit) \to H^{g+1} (A, \qit) (1))$ be the primitive cohomology. Then, we have the following:
\begin{Lem} $K_{\Theta}$ is a rational Hodge structure of level 1 and dimension 10.
\begin{proof} Since $H^{3} (\Theta)$ and $H^{3} (A)$ both have Hodge level 3, we need to show that $i^{*}: H^{3, 0} (A) \to H^{3, 0} (\Theta)$ is an isomorphism. Since this map is already injective, it will suffice to show that $h^{3, 0} (\Theta) = h^{3, 0} (A) = 4$. By adjunction, $\omega_{\Theta} \cong \mathcal{O}_{\Theta}(\Theta)$, so $h^{0} (\Theta, \mathcal{O}_{\Theta}(\Theta)) = h^{3, 0} (\Theta).$ We can use the long exact sequence to compute $h^{0} (\Theta, \mathcal{O}_{\Theta}(\Theta))$:
\begin{equation}
0 \to H^{0} (A, \mathcal{O}_{A}) \to H^{0} (A, \mathcal{O}_{A}(\Theta)) \xrightarrow{res} H^{0} (\Theta, \mathcal{O}_{\Theta}(\Theta)) \to H^{1} (A, \mathcal{O}_{A}) \to H^{1} (A, \mathcal{O}_{A}(\Theta)) = 0 \nonumber
\end{equation}
Since $\Theta$ is a principal polarization, $h^{0} (A, \mathcal{O}_{A}(\Theta)) = 1$ so that the restriction arrow is 0. Moreover, $h^{1} (A, \mathcal{O}_{A}) = 4$, so it follows that $h^{3, 0} (\Theta) = 4 = h^{3,0} (A)$. Thus, $i^{*}: H^{3, 0} (A) \to H^{3, 0} (\Theta)$ is an isomorphism and $K_{\Theta}$ has Hodge level 1. To determine the dimension of $K_{\Theta}$, we first compute $\chi(\Theta) = c_{3} (T\Theta)$. Applying the Chern polynomial to the adjunction sequence in this case, one obtains that $c_{3} (T\Theta) = - c_{1} (\mathcal{O}(\Theta))^{\cdot 4} = -4! =-24$. Using the Lefschetz hyperplane theorem, one also computes that $\chi(\Theta) = 42 - h^{3} (\Theta)$, so that $h^{3} (\Theta) = 66$. Since $h^{3} (A) = \binom{8}{3} = 56$, it follows that $K_{\Theta}$ has dimension 10. \end{proof}
\end{Lem}

Thus, $H^{*} (P, \qit)$ has Hodge level 1 when $g=4$. Now, consider the intermediate Jacobian of $K_{\Theta}$: \begin{equation} J(K_{\Theta}) = K_{\Theta, \cit}/(F^{2}K_{\Theta, \cit} \oplus K_{\Theta, \zit})\nonumber\end{equation}
This is a principally polarized Abelian variety of dimension 5, and we have an isomorphism of rational Hodge structures $H^{1} (J(K_{\Theta}), \qit) (-1) \cong H^{3} (P, \qit)$. The generalized Hodge conjecture predicts that this isomorphism arises from a correspondence $\Gamma \subset J(K_{\Theta}) \times \Theta$. The existence of $\Gamma$ was proved in \cite{IV}. One may take this a step further and ask whether $\h^{1} (J(K_{\Theta})) (-1)$ and $P$ are isomorphic as motives. Proposition \ref{3} provides a partial answer to this; i.e., we have $\h^{1} (J(K_{\Theta})) (-1) \cong P$ if $p$ acts trivially on $CH_{0} (\Theta_{L})$ for all field extensions $\cit \subset L$ (and conversely). We will need the following definition for the proof:
\begin{Def} We say that $M = (X, \pi, 0) \in \M_{k}$ has representable Chow group in codimension $i$ if there exists a smooth complete (possible reducible) curve $C$ and $\Gamma \in CH^{i} (C \times X)$ such that $CH^{i}_{alg} (M_{L}) = \pi_{L*}CH^{i}_{alg} (X_{L})$ lies in $\Gamma_{L*}CH^{1}_{alg} (C_{L}) \subset CH^{i}_{alg} (X_{L})$ for every field extension $k \subset L$.
\end{Def}

\begin{proof}[Proof of Proposition \ref{3}] Suppose that we have some Abelian variety $J$ for which $\h^{1} (J) (-1) \cong P$. Then, applying $CH^{3} (\phantom{M})$ to both sides we obtain $$p_{*}CH_{0} (\Theta) = p_{*}CH^{3} (\Theta) \cong CH^{3} (\h^{1} (J) (-1)) = CH^{2} (\h^{1} (J))$$ From \cite{DM} Theorem 2.19, we have $CH^{2} (\h^{1} (J)) = 0$ so that $p_{*}CH_{0} (\Theta) = 0$. For the converse, observe that $\Theta$ can be defined over some field $k$ which is the algebraic closure of a finitely generated over $\qit$. So, let $\Theta_{k}$ be a model for $\Theta$ over $k$. The operators used in the proof of Theorems \ref{1} and \ref{2} ($L_{\Theta}$, $\Lambda_{\Theta}$, and $\pi_{j, A}$) are well-behaved upon passage to an overfield (see \cite{DM} and \cite{K}); thus,  so is the correspondence $p_{k} \in Cor^{0} (\Theta_{k} \times \Theta_{k})$ constructed above. This means that $p_{\cit}$ coincides with $p$ (as in the statement of Proposition \ref{3}), and the assumption that $p$ acts trivially on $CH_{0}$ becomes the assumption that 
$$p_{L*}CH^{3} (\Theta_{L}) =0$$
for all overfields $k \subset L$. Now, let $P = (\Theta_{k}, p_{k}, 0)$. The task is then to find an Abelian variety $J$ over $k$ for which 
$$\h^{1} (J) (-1) \cong P$$
To this end, we begin with the following lemma:
\begin{Lem} $P$ has representable Chow group in codimension 2. 
\begin{proof}[Proof of Lemma] We will drop the subscript $k$. We use the same argument as in \cite{BS}. There is a localization sequence:
\begin{equation}\begin{CD} \displaystyle \lim_{\xrightarrow[D\subset \Theta]{}} CH^{2} (\Theta \times D) @>(id_{\Theta} \times j_{D})_{*}>> CH^{3} (\Theta \times \Theta) @>(id_{\Theta} \times K)^{*}>> CH^{3} (\Theta_{K}) @>>> 0 \label{exact}\end{CD}\end{equation} where the limit runs over all (possibly reducible) subvarieties $D$ of codimension 1 and $K = \cit(\Theta)$ is the function field of $\Theta$. We have $(id_{\Theta} \times K)^{*}\Delta_{\Theta}= \eta_{K}$, the generic point of $\Theta$. From Lemma \ref{cor-lieb} \ref{itm: h}, we have $p = p\circ \Delta_{\Theta} = (p \times id_{\Theta})_{*}\Delta_{\Theta}$ so that $$(id_{\Theta} \times K)^{*}(p) = (id_{\Theta} \times K)^{*}(p \times id_{\Theta})_{*}\Delta_{\Theta} = p_{K*}(id_{\Theta} \times K)^{*}\Delta_{\Theta}= p_{K*}(\eta_{K})$$ Since $p_{K*}(\eta_{K}) = 0$ by assumption, the exactness of (\ref{exact}) gives some subvariety $D$ and $\alpha \in CH^{2} (\Theta \times D)$ for which $p = (id_{\Theta} \times j_{D})_{*}\alpha$. After desingularizing, we can assume that $D$ is smooth (although $j_{D}$ may no longer be an inclusion). By Lemma \ref{cor-lieb} \ref{itm: a}, we have $$p =  (id_{\Theta} \times j_{D})_{*}\alpha = \Gamma_{j_{D}} \circ \alpha$$ Thus $p_{L*}CH^{2}_{alg} (\Theta_{L}) \subset j_{D*}CH^{1}_{alg} (D_{L})$. By the representability of the Picard functor, this means there is some smooth complete $C$ and some $\Gamma \in CH^{1} (C \times D)$ such that $\Gamma_{L*}CH^{1}_{alg} (C_{L}) = CH^{1}_{alg} (D_{L})$ for all field extension $k \subset L$. This proves the lemma.
\end{proof}
\end{Lem}
Thus, we see that the Chow group of $P$ is representable in every codimension. By \cite{V} Theorem 3.4, it follows that the motive of $P$ decomposes as $$\bigoplus \mathds{1}({i})^{\oplus n_{i}} \oplus \h^{1} (J_{i}) (-i)$$ for integers $n_{i}$ and Abelian varieties $J_{i}$. Since the cohomology of $P$ is 0 in all but degree 3, this means that $P \cong \h^{1} (J) (-1)$ for some Abelian variety $J$. This gives the proposition.
\end{proof}
\begin{Rem} A more refined version of Proposition \ref{3} is that the Abelian variety can be taken to be $J(K_{\Theta})$ in the above notation. Indeed, since we have $H^{1} (J(K_{\Theta}), \qit) (-1) \cong H^{3} (P, \qit) \cong H^{1} (J, \qit) (-1)$ (as rational Hodge structures), it follows that $J$ and $J(K_{\Theta})$ are isogenous.
\end{Rem}

\end{document}